\documentclass{amsproc}
\usepackage{bm, amssymb}
\newtheorem{theorem}{Theorem}[section]
\newtheorem{lemma}[theorem]{Lemma}
\newcommand\degree{\operatorname{deg}}

\DeclareMathOperator{\Mat}{M}

\theoremstyle{definition}
\newtheorem{definition}[theorem]{Definition}

\newtheorem{conjecture}[theorem]{Conjecture}

\theoremstyle{remark}

\newcommand{\mcl}[2][n]{\mathcal{#2}}

\newcommand{\defeq}{\stackrel{\text{def}}{=}}

\newcommand\cA{\mathcal{A}}

\newcommand{\bT}{\mathbf{T}}

\newcommand{\R}{\mathbb{R}}

\newcommand{\bS}{\mathbf{S}}

\def \bf {\mathbf f}

\def \bp {\mathbf p}
\def \bq {\mathbf q}

\def \bt {\mathbf t}

\def \bx {\mathbf x}  

\def \by {\mathbf y}

\def \bzero {\mathbf 0}

\def \btet {\boldsymbol{\theta}}

\numberwithin{equation}{section}

\begin{document}

\title[Diophantine approximation in positive characteristic]{The inhomogeneous Sprind\v{z}huk conjecture over a local field of positive characteristic}

\author{Arijit Ganguly}
\address{Department of Mathematics, IME Building
Indian Institute of Technology Kanpur
Kanpur, P.O.- IIT Kanpur, P.S.- Kalyanpur
District - Kanpur Nagar
Pin - 208016
Uttar Pradesh, India.}
\email{arijit.ganguly1@gmail.com}

\author{Anish Ghosh}
\address{School of Mathematics, Tata Institute of Fundamental Research, Homi Bhabha Road, Navy Nagar, Colaba, Mumbai 400005, India.}\email{ghosh@math.tifr.res.in}
\thanks{The second named author gratefully acknowledges support from a grant from the Indo-French Centre for the Promotion of Advanced Research; a Department of Science and Technology, Government of India Swarnajayanti fellowship and a MATRICS grant from the Science and Engineering Research Board.}

\date{January 1, 1994 and, in revised form, June 22, 1994.}


\begin{abstract}
We prove a strengthened version of the inhomogeneous Sprind\v{z}huk conjecture in metric Diophantine approximation, over a local field of positive characteristic. The main tool is the transference principle of Beresnevich and Velani \cite{BeVe} coupled with earlier work of the second named author \cite{G-pos} who proved the standard, i.e. homogeneous version. 
\end{abstract}
\maketitle

\section{Introduction}
The context of this paper is the metric theory of Diophantine approximation over local fields of positive characteristic. In \cite{G-pos}, the second named author proved the Sprind\v{z}huk conjectures in this setting (in fact, also in multiplicative form), here we prove the inhomogeneous variant of the conjecture. We use the inhomogeneous transference principle of Beresnevich and Velani \cite{BeVe} to transfer the homogeneous result from \cite{G-pos} and also use a positive characteristic  version of the transference principle of Bugeaud and Laurent interpolating between uniform and standard Diophantine exponents, established recently by Bugeaud and Zhang \cite{BZ}. The possibility of proving the $S$-arithmetic inhomogeneous Sprind\v{z}huk conjectures was suggested by Beresnevich and Velani (\cite{BeVe}, \S $8.4$) and the present paper realises this expectation in another natural setting, that of local fields of positive characteristic.\\

Metric Diophantine approximation on manifolds is a subject which studies the extent to which typical Diophantine properties for Lebesgue measure on $\R^n$ are inherited by smooth submanifolds or other measures. The theory began with Mahler \cite{Mah-S} who conjectured that almost every point on the Veronese curve 
is \emph{not very well approximable}. Mahler's conjecture was resolved by Sprind\v{z}huk \cite{Spr1, Spr2}, who in turn made a stronger conjecture which was resolved by Kleinbock and Margulis \cite{KM} using methods from the ergodic theory of group actions on homogeneous spaces, specifically, sharp nondivergence estimates for unipotent flows on the space of lattices. Subsequently, an $S$-arithmetic version of the conjectures were established by Kleinbock and Tomanov \cite{KT} and a positive characteristic version was established by the second named author \cite{G-pos}. Both the latter works used adaptations of the dynamical approach of Kleinbock and Margulis. In \cite{BeVe}, Beresnevich and Velani proved a transference principle which allowed them to prove an inhomogeneous versions of the Baker-Sprind\v{z}huk conjectures. We refer the reader to the above papers for more details. We will recall all the relevant concepts in the function field context in the next section.  

Following the work of Beresnevich and Velani, there have been several recent advances in inhomogeneous Diophantine approximation. In \cite{BGGV}, an inhomogeneous Khintchine type theorem was established for affine subspaces, complementing the earlier work \cite{BBV} for nondegenerate manifolds, see also \cite{GM} for more inhomogeneous results on affine subspaces. Further,  an $S$-arithmetic inhomogeneous Khintchine type theorem for nondegenerate manifolds was established by Datta and the second named author \cite{DG1}. 

\subsection{The setup}
We follow our paper \cite{GG} in setting the notation. Let $p$ be a prime and $q:= p^r$, where $r\in \mathbb{N}$, and consider the function 
field $\mathbb{F}_{q}(T)$. We define a function $|\cdot|: \mathbb{F}_{q}(T) \longrightarrow \mathbb{R}_{\geq 0}$ as follows. 
\[ |0|:= 0\,\,  \text{ and} \,\, \left|\frac{P}{Q}\right|:= e^{\displaystyle \degree P- \degree Q}
\text{ \,\,\,for all nonzero } P, Q\in \mathbb{F}_{q}[T]\,.\] 
Clearly $|\cdot|$ is a nontrivial,  non-archimedian and discrete absolute value  
in $\mathbb{F}_{q}(T)$. This absolute value gives rise to a metric on $\mathbb{F}_{q}(T)$. \\

The completion field of $\mathbb{F}_{q}(T)$ is $\mathbb{F}_{q}((T^{-1}))$, i.e. the field of Laurent series 
over $\mathbb{F}_{q}$. The absolute value of $\mathbb{F}_{q}((T^{-1}))$, which we again denote by $|\cdot |$, is given as follows. 
Let $a \in \mathbb{F}_{q}((T^{-1}))$. For $a=0$, 
define $|a|=0$. If $a \neq 0$, then we can write 
$$a=\displaystyle \sum_{k\leq k_{0}} a_k T^{k}\,\,\mbox{where}\,\,\,\,k_0 \in \mathbb{Z},\,a_k\in \mathbb{F}_{q}\,\,\mbox{and}\,\, a_{k_0}
\neq 0\,. $$
\noindent We define $k_0$ as the \textit{degree} of $a$, which will be denoted by $\degree a$,  and     
$|a|:= e^{\degree a}$. This clearly extends the absolute 
value $|\cdot|$ of $\mathbb{F}_{q}(T)$ to $\mathbb{F}_{q}((T^{-1}))$ and moreover, 
the extension remains non-archimedian and discrete. Let $\Lambda$ and $F$ 
denote $\mathbb{F}_{q}[T]$ and $\mathbb{F}_{q}((T^{-1}))$ respectively from now on. It is obvious that 
$\Lambda$ is discrete in  $F$. For any $d\in \mathbb{N}$, $F^d$ is throughout assumed to be equipped 
with the supremum norm which is defined as follows
\[||\mathbf{x}||:= \displaystyle \max_{1\leq i\leq n} |x_i|\text{ \,\,for all \,} \mathbf{x}=(x_1,x_2,...,x_d)\in F^{d}\,,\]

\noindent and with the topology induced by this norm. Clearly $\Lambda^n$ is discrete in $F^n$. Since the topology on $F^n$ considered here 
 is the usual product topology on $F^n$, it follows that  $F^n$ is locally compact as $F$ is locally compact. Let $\lambda$ be the Haar measure on $F^d$ which takes the value 1 on the closed unit ball $||\mathbf{x}||=1$.\\ 
 
Diophantine approximation in the positive characteristic setting consists of approximating elements 
in $F$ by `rational' elements, i.e. those from $\mathbb{F}_{q}(T)$. This 
subject has been extensively studied, beginning with work of E. Artin \cite{Artin} who developed the theory of continued fractions, and continuing with Mahler who developed Minkowski's geometry of numbers in 
function fields and Sprind\v{z}uk who, in addition to proving the analogue of Mahler's conjectures, 
also proved some transference principles in the function field setting (see \cite{Spr1}). The subject has also received 
considerable attention of late, we refer the reader to \cite{deM, Las1} for overviews and to \cite{AGP,  GR, GG2, Kris, KN} for a necessarily incomplete set of references. 

In the present paper we prove an inhomogeneous analogue of the Sprind\v{z}uk conjectures, our main result is an upper bound for inhomogeneous Diophantine exponents.

\begin{theorem}\label{main_th}
Let $U\subseteq F^d$ be open and $\mathbf{f}: U\longrightarrow F^n$ be a $(C,\alpha_0)-good$ map, for some $C, \alpha_0 > 0$, and assume that $(\mathbf{f},\lambda)$ is nonplanar. Then, for every $\theta \in F$, and $\lambda$ almost every $\bx \in U,$
$$\omega(\mathbf{f}(\bx), \theta) \leq 1.$$
\end{theorem} 

We also establish the corresponding lower bound.

\begin{theorem}\label{lb}
Let $U\subseteq F^d$ be open and $\mathbf{f}: U\longrightarrow F^n$ be a $(C, \alpha_0)-good$ map, for some $C, \alpha_0 > 0$, and assume that $(\mathbf{f},\lambda)$ is nonplanar. Then, for every $\theta \in F$, and $\lambda$ almost every $\bx \in U,$ 
$$\omega(\mathbf{f}(\bx), \theta) \geq 1.$$
\end{theorem}

Remarks:
\begin{enumerate}
\item Note that the exceptional set of $\bx$ for which the inequalities in Theorems \ref{main_th} and \ref{lb} need not hold depends on the inhomogeneous parameter $\theta$.\\
\item The relevant definitions are made in the next section. A main example to keep in mind is the original setup of Diophantine approximation on manifolds, i.e. if $\mathbf{f} = (f_1, \dots, f_n)$ where the $f_i$'s are analytic and $1, f_1, \dots, f_n$ are linearly independent over $F$, then $\mathbf{f}$ is $(C, \alpha)$ good for some $C, \alpha$ and nonplanar. More generally, if $\mathbf{f}$ is a smooth nondegenerate map, then it is $(C, \alpha)$-good as well as nonplanar. The notions of $(C, \alpha)$ good functions and nondegenerate maps were introduced by Kleinbock and Margulis \cite{KM}.\\
\item The homogeneous analogue of Theorem \ref{main_th} was proved in \cite{G-pos} (Theorem 3.7), the lower bound is a consequence of Dirichlet's theorem.\\ 
\item In \cite{BBV}, and subsequently in \cite{BGGV} a more general problem is considered where the inhomogeneous term is also allowed to vary. It should be possible to incorporate this improvement into Theorem \ref{main_th}.\\
\item The next five sections deal with the proof of the main theorem. Sections 2 and 3 give the necessary prerequisites, in section 4 the lower bounds for Diophantine exponents are obtained and in section 6, the corresponding upper bounds. The final section is devoted to open questions and future possibilities for research.  
\end{enumerate}

\section*{Acknowledgements} Part of this work was done when both authors were visiting ICTS Bengaluru. We thank the institute for its hospitality and excellent working conditions. We thank the referee for many helpful suggestions which have improved the paper.

\section{Homogeneous and Inhomogeneous Diophantine exponents}
The theory of Diophantine approximation in positive characteristic begins with Dirichlet's theorem, which we now recall.
\begin{theorem}\label{thm:1.1}(Theorem 2.1 \cite{GG}) 
Let $m,n \in \mathbb{N}$, $k=m+n$ 
and $$\mathfrak{a}^{+}:= \{\mathbf{t}:=(t_1,t_2,\dots,t_k)\in \mathbb{Z}_{+}^{k}\,:\,\displaystyle \sum_{i=1}^{m} t_{i} =\displaystyle \sum_{j=1}^{n}t_{m+j}\}\,.$$ 
Consider $m$ linear forms $Y_{1},Y_{2},\dots,Y_{m}$  over $F$ in $n$ variables. Then for any $\mathbf{t}\in \mathfrak{a}^+$, 
there exist solutions $\mathbf{q}=(q_1,q_2,\dots,q_n)\in \Lambda^{n} \setminus \{\mathbf{0}\}$ and $\mathbf{p}=(p_1,p_2,\dots,p_m)\in \Lambda^m$ of the 
following system of inequalities 
\begin{equation}\label{eqn:1.1}\left \{ \begin{array}{rcl} |Y_{i}\mathbf{q}-p_i|\textless e^{-t_{i}} &\mbox{for} &i=1,2,\dots,m \\ |q_j|\leq e^{t_{m+j}}&\mbox{for} &j=1,2,\dots,n\,. \end{array} \right.  
\end{equation}
\end{theorem}

We will consider only unweighted Diophantine approximation in this paper, so $t_1= \dots = t_m = 1/m$ and $t_{n+1} = \dots = t_{n} = 1/n$. We denote by $\Mat_{m \times n} (F)$, the vector space of $m \times n$ matrices with entries from $F$ equipped with the supremum norm. In view of Theorem \ref{thm:1.1}, it is natural to define exponents of Diophantine approximation as follows. Let $X \in \Mat_{m\times n}(F)$ and $\btet\in F^m$. The \emph{inhomogeneous exponent}, $\omega(X, \btet)$ of $X$, is the supremum of the real numbers $\omega$ for which, for arbitrarily large  $T\in \mathbb{N}$, the inequalities
\begin{equation}\label{equ-def-w}
   \|X\bq-\bp-\btet\| <e^{-\frac{n}{m}\omega T},  \qquad  \|\bq\|< e^{T},
\end{equation}
have a  solution $(\bp,\bq)\in \Lambda^m\times(\Lambda^n\setminus \{\bzero\})$.
The \emph{uniform inhomogeneous exponent}, $\hat{\omega}(X,\btet)$, is the supremum of the real numbers $\hat{\omega}$ for which, for all sufficiently $T\in \mathbb{N}$, the inequalities
\begin{equation*}
   \|X\bq-\bp-\btet\|<e^{-\frac{n}{m}\hat{\omega} T}, \qquad  \|\bq\|<e^T, 
\end{equation*}
have a solution $(\bp,\bq)\in \Lambda^m\times(\Lambda^n\setminus \{\bzero\})$. \\

In this paper, we will adopt the point of view of Diophantine approximation of single linear forms, i.e. we will assume that $\by \in F^n$ where $F^n$ is identified with $\Mat_{1 \times n }(F)$  as opposed to simultaneous Diophantine approximation where one considers $\by \in \Mat_{n \times 1}$.\\

If $\theta = 0$, then the corresponding Diophantine exponent $\omega(\mathbf{y}) := \omega(\mathbf{y}, 0)$  (resp. $\hat{\omega}(\by)$) is called the homogeneous Diophantine exponent. By Dirichlet's theorem stated above, $\omega(\mathbf{y}) \geq 1$ for every $\mathbf{y} \in F^{n}$. We are following the normalisation in \cite{BeVe} rather than the one used in \cite{KM, G-pos} according to which the critical exponent is $n$.\\

The Borel-Cantelli lemma implies that $\omega(\mathbf{y}) = 1$ for $\lambda$ almost every $\mathbf{y} \in F^{n}$. It is therefore natural to define $\mathbf{y} \in F^{n}$ to be \emph{very well approximable} if $\omega(\mathbf{y}) > 1$. Sprind\v{z}huk \cite{Spr2} proved that for $\lambda$ a.e. $x \in F$,
\begin{equation}\label{veronese}
\mathbf{f}(x) := (x, x^2, \dots, x^n)
\end{equation}
is not very well approximable, thereby settling the positive characteristic analogue of Mahler's conjecture. A special case of the theorems proved in this paper is that for  every $\theta \in F$,  $\omega(\mathbf{f}(x), \theta) = 1$ for almost every $x$. Following \cite{BeVe} we may define inhomogeneously extremal measures as follows.

\begin{definition}
Let $\mu$ be a measure supported on a subset of $F^{n}$. We say that $\mu$ is inhomogeneously extremal if for all $\theta\in F$,
$$\omega(\mathbf{y}, \theta) = 1 \text{ for } \mu\text{ a.e. }\mathbf{y} \in F^{n}.$$
\end{definition}

Then our main theorems can be restated as follows:

\begin{theorem}
Let $U\subseteq F^d$ be open and $\mathbf{f}: U\longrightarrow F^n$ be a $(C,\alpha_0)$-good map, for some $C,\alpha_0 > 0$, and assume that $\mathbf{f}$ is nonplanar. Then $\mathbf{f}_{*}\lambda$ is inhomogeneously extremal.
\end{theorem}

\section{Good and nonplanar maps}
We recall the following definitions and results from \cite[\S 1 and 2]{KT}. For the sake of generality, we assume $X$ is a Besicovitch metric space, $U\subseteq X$ is  open, $\nu$ is a Radon measure
on $X$, $(\mcl{F},|\cdot|)$ is a valued 
field and $f: X\longrightarrow\mcl{F}$ is a given function such that $|f|$ is measurable. Recall that a metric space $X$ is called \emph{Besicovitch} \cite{KT} if there exists a constant $N_X$ such that the following holds: for any bounded subset $A$ of $X$ and for any family $\mathcal{B}$ of nonempty open balls in $X$ such that
$$ \forall x \in A \text{ is a center of some ball of } B,$$
 there is a finite or countable subfamily $\{B_i\}$ of $B$ with
 $$ 1_A \leq \sum_{i}1_{B_i} \leq N_X. $$
 For any $B\subseteq X$, we set
$$ ||f||_{\nu,B} := \displaystyle \sup_{x\in B\cap \text{ supp }(\nu)} |f(x)|.$$

\begin{definition}\label{defn:C,alpha}
 For $C, \alpha \textgreater\,0$, $f$ is said to be $(C,\alpha)$-good on $U$ with respect to $\nu$ if for every ball $B\subseteq U$
 with center in $\text{supp }(\nu)$, one has
 \[\nu(\{x\in B: |f(x)|\,\textless \varepsilon\})\leq C\left(\frac{\varepsilon}{||f||_{\nu,B}}\right)^{\alpha} \nu(B)\,.\]
\end{definition}
The following properties are immediate from Definition \ref{defn:C,alpha}. 
\begin{lemma} \label{lem:C,alpha} Let $X,U,\nu, \mcl{F}, f, C,\alpha,$ be as given above. Then one has 
\begin{enumerate}
   \item $f$ is  $(C,\alpha)$-good on $U$ with respect to $\nu \Longleftrightarrow$ so is $|f|$.
   \item $f$ is $(C,\alpha)$-good on $U$  with respect to $\nu \Longrightarrow$ so is $cf$ for all $c \in \mcl{F}$.
   \item \label{item:sup} $\forall~i\in I, f_i$ are $(C,\alpha)$-good on $U$ with respect to $\nu$ and $\sup_{i\in I} |f_i|$ is measurable $\Longrightarrow$ so is $\sup_{i\in I} |f_i|$.
   \item $f$ is $(C,\alpha)$-good on $U$ with respect to $\nu$ and $g :V \to \mathbb{R}$ is a continuous function such 
   that $c_1\leq |\frac{f}{g}|\leq c_2$ 
 for some $c_1,c_2 \,\textgreater \,0\Longrightarrow g$ is $(C(\frac{c_2}{c_1})^{\alpha},\alpha)$-good on $U$ with respect to $\nu$.
 \item Let $C_2 \,\textgreater \,1$ and 
 $\alpha _2\,\textgreater\,0$. Then $f$ is $(C_1,\alpha_1)$-good on $U$ with respect to $\nu$ and $C_1 \leq C_2, \alpha_2 \leq \alpha_1 \Longrightarrow f$ is
 $(C_2,\alpha_2)-good$ on $V$ with respect to $\nu$.
  \end{enumerate}
\end{lemma}  
We say a  map $\mathbf{f}=(f_1, f_2, \dots, f_n)$ from $U$ to $\mcl{F}^n$, where $n\in \mathbb{N}$,
is $(C,\alpha)$-good on $U$ with respect to $\nu$, or simply $(\mathbf{f},\nu)$ is 
$(C,\alpha)$-good on $U$,  if every $\mcl{F}$-linear combination of $1,f_1, \dots, f_n$ 
is $(C,\alpha)$-good on $U$ with respect to $\nu$. 

\begin{definition} Let $\mathbf{f}=(f_1, f_2,\dots, f_n)$ be a map from $U$ to $\mcl{F}^n$, where $n\in \mathbb{N}$. We say that $(\mathbf{f},\nu)$ is \emph{nonplanar} at a given point $x_0\in U$ if for any ball $B$ with centered at $x_0$, the restrictions of the functions  $1,f_1,\dots, f_n$ 
on $B\cap \text{ supp }(\nu)$ are linearly independent over $\mathcal{F}$. If $(\mathbf{f}, \nu)$ is nonplanar at $\nu$ almost every point of $U$, then it is called $\emph{nonplanar}$. We also simply say $\mathbf{f}$ is \emph{nonplanar} when there is no possibility of confusion. 
 \end{definition}
A typical example is provided by $\mathbf{f}=(f_1, f_2,\dots, f_n)$ where $1, f_1, \dots, f_n$ are smooth and linearly independent on $U$. Such a map has been called \emph{nondegenerate} by Kleinbock and Margulis. \\

For $m\in \mathbb{N}$ and a ball $B=B(x;r)\subseteq X$, where $x\in X$ and $r\,\textgreater\,0$, we shall use the notation
$3^mB$ to denote the ball $B(x;3^mr)$. Finally, we will need the notion of a doubling measure.
\begin{definition}
  The measure $\nu$ is said to be \emph{doubling} on $U$ if there exists $D>0$ such that for every ball $B$ with center 
 in $\text{supp }(\nu)$ such that $2B\subseteq U$, one has \[\frac{\nu(2B)}{\nu(B)}\leq D\,.\]
\end{definition}

\section{Transference principles and lower bounds}
The lower bound will follow immediately from two Diophantine transference principles. The following result was proved by Bugeaud and Zhang \cite{BZ} and constitutes a positive characteristic version of the transference principle of Bugeaud and Laurent \cite{BL}.
\begin{theorem}(Theorem 1.2, \cite{BZ})
Let $X \in F^{m \times n}$. Then for all $\btet \in F^m$, we have
\begin{equation}\label{BL}
\omega(X, \btet) \geq \frac{1}{\hat{\omega}(X^{t})} \text{ and } \hat{\omega}(X, \btet) \geq \frac{1}{\omega(X^t)}.
\end{equation}
with equalities for almost every $\btet$.
\end{theorem}
 
 We will also need a positive characteristic version of Dyson's transference principle  \cite{Dyson} which can be formulated as follows.
 \begin{theorem}
 For $\mathbf{y} \in F^n,$
 $$\omega(\mathbf{y}) = 1 \text{ if and only if } \omega(^{t}\mathbf{y}) = 1.$$
 \end{theorem}

We omit the short proof which can be obtained by a verbatim repetition of the proof in \cite{Dyson}, or the more recent, more general version proved in Theorem 1.7 in  \cite{CGGMS}.\\

It is now easy to complete the proof of the lower bound Theorem \ref{lb}.

\begin{proof}
Under the hypothesis of Theorem \ref{lb}, using Theorem $3.7
$ of \cite{G-pos}, we have that for $\lambda$ almost every $\bx$, $\omega(\mathbf{f}(\bx)) = 1$. Set $y = \mathbf{f}(\bx)$, then by Dyson's transference principle, $\omega(^{t}\mathbf{y}) = 1$. By Dirichlet's theorem, $\hat{\omega}(\by) \geq 1$ and the trivial inequality $$\omega(\by, \theta) \geq \hat{\omega}(\by, \theta) \geq 0,$$ 
applied to $^{t}\mathbf{y}$ and $\theta = 0$ we get that $\omega(^{t}\mathbf{y}) = 1$. Finally, by (\ref{BL}), we get that $\omega(\by, \theta) \geq 1$ which completes the proof.

\end{proof}

\section{The Transference principle of Beresnevich-Velani}\label{transference}
In this section we state the  inhomogeneous transference principle of Beresnevich and Velani from \cite[Section 5]{BeVe} which will allow 
us to convert our inhomogeneous problem to the homogeneous one. Let $(\Omega, d)$ be a locally compact metric space. Given two countable indexing sets $\mathcal{A}$ and $\mathbf{T}$, let H and I be two maps from $\bT \times \cA \times \R_{+}$ into the set of open subsets of $\Omega$ such that

\begin{equation}\label{H_fn}
 H~:~(t, \alpha, \varepsilon) \in \bT \times \cA \times \R_{+} \to H_{\mathbf{t}}(\alpha, \varepsilon) 
\end{equation}
\\

and
\begin{equation}\label{I_fn}
 I~:~ (t, \alpha, \varepsilon) \in \bT \times \cA \times \R_{+} \to I_{\mathbf{t}}(\alpha, \varepsilon). 
\end{equation}  
\\
Furthermore, let
\begin{equation}\label{defH}
 H_{\bt} (\varepsilon) := \bigcup_{\alpha \in \cA} H_{\mathbf{t}}(\alpha, \varepsilon) \text{ and }  I_{\bt} (\varepsilon) := \bigcup_{\alpha \in \cA} I_{\mathbf{t}}(\alpha, \varepsilon).
\end{equation}

Let $\Psi$ denote a set of functions $\psi: \bT \to \R_{+}~:~\bt \to \psi_{\bt}$. For $\psi \in \Psi$, consider the limsup sets

\begin{equation}\label{deflambda}
\Lambda_{H}(\psi) = \limsup_{\bt \in \bT} H_{\bt}(\psi_{\bt}) \text{ and } \Lambda_{I}(\psi) = \limsup_{\bt \in \bT} I_{\bt}(\psi_{\bt}).
\end{equation}

The sets associated with the map $H$ will be called homogeneous sets and those associated with the map $I$, inhomogeneous sets. We now come to two important properties connecting these notions.

\subsection*{The intersection property} The triple $(H, I, \Psi)$ is said to satisfy the intersection property if, for any $\psi \in \Psi$, there exists $\psi^{*} \in \Psi$ such that, for all but finitely many $\bt \in \bT$ and all distinct $\alpha$ and $\alpha'$ in $\cA$, we have that
\begin{equation}\label{inter}
I_{\bt}(\alpha, \psi_{\bt}) \cap  I_{\bt}(\alpha', \psi_{\bt}) \subset H_{\bt}(\psi^{*}_{\bt}).
\end{equation}

\subsection*{The contraction property}  Let $\mu$ be a finite, non atomic, doubling measure supported on a bounded subset $\mathbf{S}$ of $\Omega$. 
We say that $\mu$ is contracting with respect to $(I, \Psi)$ if, for any
$\psi \in \Psi$, there exists $\psi^{+}\in \Psi$ and a sequence of positive numbers $\{k_{\bt}\}_{\bt \in \bT}$ satisfying
\begin{equation}\label{conv}
\sum_{\bt \in \bT}k_{\bt} < \infty,
\end{equation}
such that, for all but finitely $\bt \in \bT$ and all $\alpha \in \cA$, there exists a collection $C_{\bt, \alpha}$ of balls $B$
centred at $\mathbf{S}$ satisfying the following conditions:
\begin{equation}\label{inter1}
\bS \cap I_{\bt}(\alpha, \psi_{\bt}) \subset \bigcup_{B \in C_{\bt, \alpha}} B
\end{equation}

\begin{equation}\label{inter2}
\bS \cap \bigcup_{B \in C_{\bt, \alpha}} B \subset  I_{\bt}(\alpha, \psi^{+}_{\bt})
\end{equation}
and

\begin{equation}\label{inter3}
\mu(5B \cap  I_{\bt}(\alpha, \psi_{\bt})) \leq k_{\bt} \mu(5B).
\end{equation}

We are now in a position to state Theorem $5$ from \cite{BeVe}. 
\begin{theorem}\label{transfer}
Suppose that $(H, I, \Psi)$ satisfies the
intersection property and that $\mu$ is contracting with respect to $(I, \Psi)$. Then
\begin{equation}\label{eq:transfer1}
\mu(\Lambda_{H}(\psi))=0 ~\forall~\psi \in \Psi  \Rightarrow \mu(\Lambda_{I}(\psi)) = 0 ~\forall~\psi \in \Psi.
\end{equation}
\end{theorem}

\section{Proof of Theorem \ref{main_th}}
Fix $\theta\in F$. It is enough to show that for any open ball $V\subseteq U$ such that $5V\subseteq U$, 
$
 \omega(\mathbf{f}(\mathbf{x}),\theta)\leq 1\text{ for }\lambda\text{ almost all } \mathbf{x}\in V.
$ In fact, we prove 
\[\forall ~\omega>1, \lambda(\{\mathbf{x}\in V: \omega(\mathbf{f}(\mathbf{x}),\theta)> \omega\})=0.
\]

For each $(t,\alpha=(p,\mathbf{q}),\varepsilon)\in \mathbb{N}\times (\Lambda\times\Lambda^n\setminus \{0\})\times \mathbb{R}_{+}$, we set 
\[ I_{\mathbf{t}}(\alpha, \varepsilon)\defeq \{\mathbf{x}\in V: |\mathbf{f}(\mathbf{x})\cdot\mathbf{q}+p+\theta|<\varepsilon, ||\mathbf{q}||\leq e^t\}, \] and 
\[ H_{\mathbf{t}}(\alpha, \varepsilon)\defeq \{\mathbf{x}\in V: |\mathbf{f}(\mathbf{x})\cdot\mathbf{q}+p|<\varepsilon, ||\mathbf{q}||\leq e^t\}.\] 
Let $\Psi$ denote the collection of functions $\psi_{\omega}: \mathbb{N}\rightarrow \mathbb{R}, t\mapsto \frac{1}{e^{n\omega t}}$, for $\omega >1$. We denote the restriction of $\lambda$ to $V$ by $\mu$  and thus it is supported on $V$. \\

Since 
$\forall ~\omega>1, \{\mathbf{x}\in V: \omega(\mathbf{f}(\mathbf{x}),\theta)> \omega\})\subseteq \Lambda_I(\psi_{\omega})
$ so, it suffices to show that $\lambda(\Lambda_I(\psi_{\omega}))=0$ for any $\omega>1$. 
Theorem $3.7$ in \cite{G-pos} implies that 
$$\forall ~\omega >1,~\lambda(\Lambda_{H}(\psi_{\omega}))=0.$$ Therefore to prove Theorem \ref{main_th}, in view of the Theorem \ref{transfer}, we only need to verify the intersection and contraction properties. These will be performed in the following two subsections. 

\subsection{Verification of the intersection property}
Let $t\in \mathbb{N}, \alpha=(p,\mathbf{q}),\alpha'=(p',\mathbf{q'})\in \Lambda \times \Lambda^{n}\setminus \{0\}$ with $\alpha\neq \alpha'$ and $\omega>1$. If at least one of $||\mathbf{q}||$ and  $||\mathbf{q'}||$ is $>e^t$, then there is nothing to prove. Otherwise, the ultrametric property yields that if $\mathbf{x}\in I_t(\alpha,\psi_{\omega}(t))\cap I_t(\alpha',\psi_{\omega}(t))$ then
\begin{equation}\label{inter}
 \begin{array}{rcl}
|\mathbf{f}(\mathbf{x})\cdot (\mathbf{q}-\mathbf{q'})+(p-p')|\leq \max\{|\mathbf{f}(\mathbf{x})\cdot \mathbf{q}+p+\theta|, |\mathbf{f}(\mathbf{x})\cdot \mathbf{q'}+p'+\theta|\}\leq  \frac{1}{e^{n\omega t}}.
\end{array}
\end{equation} Note that if $\mathbf{q}=\mathbf{q}'$, then $|p-p'|\leq \frac{1}{e^{n\omega t}}$ and so $p=p'$ which is impossible. Hence, it follows from (\ref{inter}) that $I_t(\alpha,\psi_{\omega}(t))\cap I_t(\alpha',\psi_{\omega}(t))\subseteq H_t(\alpha-\alpha', \psi_{\omega}(t))$. 

\subsection{Verification of the contraction property}
Fix $\alpha \in \Lambda \times \Lambda^n\setminus \{0\}$. We observe that, for any $t\in\mathbb{N}$, 
$I_t(\alpha,\psi_{\omega}(t))\subseteq I_t(\alpha,\psi_{\frac{\omega+1}{2}}(t))$ and 
\begin{equation}\label{f_good}
 \mu(I_t(\alpha,\psi_{\frac{\omega+1}{2}}(t)))\leq \mu\left(\{\mathbf{x}\in V: |\mathbf{f}(\mathbf{x})\cdot \mathbf{q}+p+\theta|<\frac{1}{e^{\frac{\omega+1}{2}nt}}\}\right)\ll\frac{1}{e^{\frac{\omega+1}{2}nt\alpha_0}}\mu(V), 
\end{equation}
 since $\mathbf{f}$ is $(C,\alpha_0)$-good on $U$. From the nonplanarity of $\mathbf{f}$, we have $$\displaystyle \inf_{\alpha} \sup_{\mathbf{x}\in U} |\mathbf{f}(\mathbf{x})\cdot \mathbf{q}+p+\theta|>0. $$ So the absolute constant appearing in the last inequality of (\ref{f_good}) can be made independent of $\alpha$. Thus it turns out from (\ref{f_good}) that, for all sufficiently large $t$, 
 
 \begin{equation}\label{subset}
  I_t(\alpha,\psi_{\frac{\omega+1}{2}}(t))\subsetneqq V \text{ for all } \alpha. 
 \end{equation}
  
 For any $t$ that satisfies (\ref{subset}) and all $\alpha$, we now construct a collection of balls $C_{t,\alpha}$ centered in $V$ which makes (\ref{inter1})-(\ref{inter3}) hold. If $I_t(\alpha,\psi_{\omega}(t))=\emptyset$ then we set $C_{t,\alpha}$ as the empty collection and consequently, (\ref{inter1})-(\ref{inter3}) become trivial. Suppose $I_t(\alpha,\psi_{\omega}(t))$ is nonempty. Let $\mathbf{x}\in I_t(\alpha,\psi_{\omega}(t))$. Since $I_t(\alpha,\psi_{\frac{\omega+1}{2}}(t))$ is open, there exists a ball $B'(\mathbf{x})$ with center $\mathbf{x}$ such that $B'(\mathbf{x})\subseteq I_t(\alpha,\psi_{\frac{\omega+1}{2}}(t))$. We can scale it and denote it by $B(\mathbf{x})$, due to (\ref{subset}), in such a way that 
 \begin{equation}\label{construction}
  B(\mathbf{x})\subseteq I_t(\alpha,\psi_{\frac{\omega+1}{2}}(t))\nsupseteq V\cap 5B(\mathbf{x}).
 \end{equation} It is also clear from the construction that $5B(\mathbf{x})\subseteq 5V\subseteq U$.
Consider 
$$C_{t,\alpha}\defeq \{B(\mathbf{x}): \mathbf{x}\in I_t(\alpha,\psi_{\omega}(t))\}.$$ 
The conditions (\ref{inter1}) and (\ref{inter2}) are obvious. \\

Define $F_{\alpha}:U\longrightarrow \mathbb{R}, F_{\alpha}(\mathbf{x})=|\mathbf{f}(\mathbf{x})\cdot\mathbf{q}+p+\theta|,~\forall ~ \mathbf{x}\in U$ and let $B\in C_{t,\alpha}$. 
 By the last inequality given in (\ref{construction}), we see that 
\begin{equation}\label{above}
 \sup_{\mathbf{x}\in 5B}F_{\alpha}(\mathbf{x})\geq \sup_{\mathbf{x}\in 5B\cap V}F_{\alpha}(\mathbf{x})\geq \frac{1}{e^{\frac{\omega+1}{2}nt}}.
\end{equation} Furthermore, one has
\begin{equation}\label{estimate}
 \sup_{\mathbf{x}\in 5B\cap I_t(\alpha,\psi_\omega (t))}F_{\alpha}(\mathbf{x})
<\frac{1}{e^{n\omega t}}\leq \frac{1}{e^{n\omega t}}\times e^{\frac{\omega+1}{2}nt}\sup_{\mathbf{x}\in 5B}F_{\alpha}(\mathbf{x})=\frac{1}{e^{\frac{w-1}{2}nt}}\sup_{\mathbf{x}\in 5B}F_{\alpha}(\mathbf{x}),
\end{equation} due to (\ref{above}).
 Hence, from (\ref{estimate}) and the assumption that $\mathbf{f}$ is $(C,\alpha_0)$-good on $U$, it follows now that
\begin{equation}\label{int3}
 \begin{array}{rcl}
  \mu(5B\cap I_t(\alpha,\psi_\omega (t)))=\lambda(5B\cap I_t(\alpha,\psi_\omega (t)))\\ \leq \lambda\left(\displaystyle\{\mathbf{x}\in 5B: F_{\alpha}(\mathbf{x})<\frac{1}{e^{\frac{w-1}{2}nt}}\sup_{\mathbf{x}\in 5B}F_{\alpha}(\mathbf{x})\}\right)\\ \displaystyle \leq \frac{C}{e^{\frac{w-1}{2}nt\alpha_0}}\lambda(5B).
 \end{array}
\end{equation}
Since $5B\cap V= V \text{ or }5B$, accordingly as $V\subseteq 5B$ or $5B\subseteq V$, so we have  $\mu(5B)=\lambda(V)\text{ or }\lambda(5B)$. In the first case, we obtain $\lambda(5B)\leq  \lambda(5V)=5^n\lambda(V)=5^n\mu(5B)$, and $\lambda(5B)=\mu(5B)$ in the later. Thus in either case, we see that  $\lambda(5B)\leq 5^n \mu(5B)$. In view of this and (\ref{int3}), the condition (\ref{inter3}) of the contraction property is obvious as soon as we set 
\[k_t\defeq \frac{5^n C}{e^{\frac{w-1}{2}nt\alpha_0}},~\forall ~ t\gg1.\]

\section{Further directions}

In this section, we mention some directions for future research. 
\subsection{One vs almost every dichotomies}
In \cite{Kleinbock-dichotomy}, D. Kleinbock proved a remarkable \emph{dichotomy} for Diophantine exponents. A special case of his results implies that if a connected analytic manifold $\mathcal{M} \subset \mathbb{R}^n$ has one not very well approximable point, then almost every point on $\mathcal{M}$ is not very well approximable. In \cite{DG2}, a $p$-adic version of this result was obtained. It is natural to ask if inhomogeneous analogues of Kleinbock's results hold. In other words, we propose

\begin{conjecture}
Let $\mathcal{M} \subset \mathbb{R}^n$ be a connected analytic manifold. Suppose there exists $\mathbf{x} \in \mathcal{M}$ such that for every $\theta \in \mathbb{R}$,
\begin{equation}\label{conj1}
\omega(\bx, \theta) = 1.
\end{equation}
Then $\mathcal{M}$ is inhomogeneously extremal.
\end{conjecture}
This conjecture can of course be formulated over any local field as well as in the multiplicative setting. It should be noted that Kleinbock's technique does not seem to be directly applicable in the inhomogeneous setting. 

\subsection{Diophantine approximation on limit sets}
Beginning with pioneering work of Patterson \cite{Pat}, the theory of metric Diophantine approximation in the context of dense orbits of geometrically finite Kleinian groups has developed into a full fledged theory. Recently, in \cite{BGSV}, a theory of metric Diophantine approximation on manifolds was developed in the context of Kleinian groups. Namely, questions of inheritance of Diophantine properties for proper subsets of the limit set of a Kleinian group were investigated. This theory has a natural counterpart in positive characteristic; where one considers orbits of discrete subgroups of $G(k)$ for algebraic groups $G$ defined over $k$ on the boundary of the Bruhat-Tits building. It would be interesting to obtain a ``manifold" theory in this context analogous to \cite{BGSV}.  

\subsection{Friendly and nonplanar measures and multiplicative Diophantine approximation}
It should be possible to extend our main Theorem to a wider class of measures, namely strongly contracting measures as considered by Beresnevich and Velani \cite{BeVe}. This class of measures includes friendly measures as defined by Kleinbock, Lindenstrauss and Weiss \cite{KLW}. Though we do not discuss this here, in fact Theorems \ref{main_th} and \ref{lb} should hold for a wider class of measures, the so called \emph{strongly contracting measures} as introduced by Beresnevich and Velani, a category which includes the important class of \emph{friendly} measures introduced earlier by Kleinbock, Lindenstrauss and Weiss \cite{KLW}. It should also be possible to extend the main Theorem to the setting of multiplicative Diophantine approximation, thereby obtaining an inhomogeneous analogue of Baker's strong extremality conjecture.

\subsection{Khintchine-Groshev type theorems}
In \cite{Kris}, S. Kristensen proves an asymptotic formula for the number of solutions to inhomogeneous Khintchine type inequalities for matrices with entries in $F$, thereby obtaining an analogue of W. Schmidt's results \cite{S1, S2} in the positive characteristic setting. While this generality seems out of reach at present in the context of manifolds, it would be interesting to prove a qualitative result, namely homogeneous and inhomogeneous Khintchine type theorems for smooth manifolds in the positive characteristic setting. These would constitute function field analogues of the work of Bernik, Kleinbock and Margulis \cite{BKM} who proved the convergence Khintchine theorem for smooth nondegenerate manifolds, and Beresnevich, Bernik, Kleinbock and Margulis \cite{BBKM} who proved the divergence case. In the inhomogeneous case, the convergence and divergence khintchine type theorems were proved by Badziahin, Beresnevich and Velani \cite{BBV}.

\subsection{Affine subspaces and their nondegenerate submanifolds}
The results in the present paper have to do with nondegenerate manifolds. At the other end of the spectrum lie affine subspaces, the study of whose Diophantine properties involves subtle considerations concerning the slope of the subspace. There has been considerable work in this area recently, cf.  \cite{Kleinbock-extremal, Kleinbock-exponent, G1, G-thesis, G-div, G-mult, G-monat}. We refer the reader to \cite{G-handbook} for a survey of this subject.  It would be interesting to obtain function field analogues of these results, both homogeneous and inhomogeneous.

\bibliographystyle{amsalpha}

\end{document}